\documentstyle{amsppt}
\magnification=1200
\hsize=6.5truein
\vsize=8.9truein
\topmatter
\title Bloch's Conjecture and Chow Motives
\endtitle
\author   Morihiko Saito
\endauthor
\affil RIMS Kyoto University, Kyoto 606-8502 Japan \endaffil
\endtopmatter
\tolerance=1000
\baselineskip=12pt 
\def\scirc{\raise.2ex\hbox{${\scriptstyle\circ}$}}
\def\sbull{\raise.3ex\hbox{${\scriptscriptstyle\bullet}$}}
\def\CH{\hbox{{\rm CH}}}
\def\Ext{\hbox{{\rm Ext}}}
\def\Ker{\hbox{{\rm Ker}}}
\def\Hom{\hbox{{\rm Hom}}}
\def\Alb{\text{{\rm Alb}}}
\def\HS{\text{{\rm HS}}}
\def\simto{\buildrel\sim\over\to}
\def\SameAuthor{\vrule height3pt depth-2.5pt width1cm}

\document
\noindent
Let
$ X $ be a connected smooth projective complex surface.
J. Murre [7] constructed a decomposition of Chow motives for
$ X $,
i.e.
there exist mutually orthogonal idempotents
$ \pi _{i} \in  \CH^{2}(X\times X)_{{\Bbb Q}} $ as correspondences for
$ 0 \le  i \le  4 $ such that
$ \sum _{i} \pi _{i} $ is equal to the diagonal and the action of
$ \pi _{i} $ on
$ H^{j}(X,{\Bbb Q}) $ is the identity for
$ i = j $,
and vanishes otherwise.
The decomposition is not uniquely characterized by the above properties.
See also [8].
The theory of Chow motives would be rather complicated if the following condition
is not satisfied:

\medskip
\noindent
(0.1) An idempotent of
$ \CH^{2}(X\times X)_{{\Bbb Q}} $ is zero if so is its cohomology class.

\medskip

Here we can also consider a stronger condition:

\medskip
\noindent
(0.2) $ \Ker(\CH^{2}(X\times X)_{{\Bbb Q}} \rightarrow  H^{4}(X\times X,
{\Bbb Q})(2)) $ is a nilpotent ideal.

\medskip

This is related to Beilinson's conjectures [2].
If (0.2) is true, the uniqueness of the projectors modulo inner 
automorphisms can be proved due to Beilinson.
See [7, 7.3].

Let
$ h^{i}(X) $ denote the `image' of the projector
$ \pi _{i} $ in the motivic sense.
Then
$ h^{2}(X) $ carries the Albanese kernel and the Neron-Severi group as well
as the transcendental cycles.
This is compatible with a conjecture of S.~Bloch [3] that 
the vanishing of the transcendental cycles (i.e. that of
$ p_{g})  $ would be equivalent to:

\medskip
\noindent
(0.3) The Albanese map is injective.

\medskip

Note that the noninjectivity of the Albanese map in the case
$ p_{g} \ne  0 $ is a theorem of D.~Mumford [6], and Bloch's conjecture is proved
at least if
$ X $ is not of general type [4].
See also [1], [13] for some more cases.

In [7] the second projector
$ \pi _{2} $ is actually defined as the difference between the diagonal
and the sum of the other projectors
$ \pi _{i} $.
In the case
$ p_{g} = 0 $,
we can construct explicitly a projector
$ \tilde{\pi}_{2} $ which is homologically equivalent to
$ \pi _{2} $ and is orthogonal to the other projectors
$ \pi _{i} $ 
$ (i \ne  2) $.
Then we can prove Bloch's conjecture if
$ \tilde{\pi}_{2} $ coincides with
$ \pi _{2} $ (modulo rational equivalence).
So the conjecture is reduced to (0.1).
Actually, we can show

\proclaim{{\bf 0.4.~Theorem}}
In the case
$ p_{g} = 0,  $ the above three conditions (0.1--3) are 
all equivalent, and the cube of the ideal in (0.2) is zero if it is nilpotent.
\endproclaim

The proof of (0.3) 
$ \Rightarrow  (0.2)  $ uses an argument similar to [3] together with the
bijectivity of the cycle map in the divisor case.
Combined with the above mentioned result of [4], it implies

\proclaim{{\bf 0.5.~Theorem}}
If 
$ p_{g} = 0 $ and
$ X $ is not of general type (or, if
$ X $ is as in [1], [13]), then the cube of the ideal in (0.2) is zero.
\endproclaim

We can show that the square of the ideal in (0.2) does not vanish if the 
irregularity 
$ q \,(= \dim \Gamma(X, \Omega_{X}^{1})) $ is nonzero.

In Sect. 1, we review Murre's decomposition of Chow motives, and prove 
(0.1) 
$ \Rightarrow  $ (0.3).
In Sect. 2 we show (0.3) 
$ \Rightarrow  $ (0.2) using a variant of the construction of [3].

Part of this work was done during my stay at the university of Leiden.
I would like to thank J.P. Murre for useful discussions on Chow motives that
have originated this work.
I thank also the staff of the institute for the hospitality.

\bigskip\bigskip\centerline{{\bf 1. Chow motives}}

\bigskip
\noindent
{\bf 1.1.}
{\it Correspondences.} For smooth proper complex algebraic varieties
$ X, Y $ such that
$ X $ is purely
$ n $-dimensional, we define the group of correspondences with 
rational coefficients by
$$
C^{i}(X,Y)_{{\Bbb Q}} = \CH^{n+i}(X\times Y)_{{\Bbb Q}}.
$$
For
$ \xi  \in  C^{i}(X,Y)_{{\Bbb Q}} $ and
$ \eta  \in  C^{j}(Y,Z)_{{\Bbb Q}} $,
the composition is denoted by
$ \eta \scirc\xi  \in  C^{i+j}(X,Z)_{{\Bbb Q}} $.
For
$ \zeta  \in  \CH^{i}(X)_{{\Bbb Q}} $,
let
$$
\Gamma _{\zeta } \in  C^{i}(pt,X)_{{\Bbb Q}} \,(= \CH^{i}(X)_{{\Bbb Q}})
\leqno(1.1.1)
$$
be the element defined by
$ \zeta  $.
For a morphism
$ f : X \rightarrow  Y $,
we denote by
$ \Gamma _{f} $ the graph of
$ f $ which belongs to
$ C^{0}(Y,X) $.
Sometimes we will use the notation
$$
f^{*} = \Gamma _{f},\quad f_{*} = {}^{t}\Gamma _{f}.
\leqno(1.1.2)
$$

Assume
$ X, Y $ connected.
By Hodge theory together with the K\"unneth decomposition and the 
duality, we have a canonical isomorphism
$$
{{\CH}^{1}(X\times Y{)}_{{\Bbb Q}} \over 
{pr}_{1}^{*}{\CH}^{1}(X{)}_{{\Bbb Q}}+{pr}_{2}^{*}{\CH}^{1}(Y{)}_{{\Bbb Q}}} = 
\Hom_{\HS}(H^{2n-1}(X,{\Bbb Q})(n-1),H^{1}(Y,{\Bbb Q})),
\leqno(1.1.3)
$$
where the right-hand side is the group of morphisms of Hodge structures.
See also [7], [12].
Let
$ \xi  \in  \CH_{0}(X)_{{\Bbb Q}}, \xi ' \in  \CH_{0}(Y)_{{\Bbb Q}}
 $ with degree one.
Then the left-hand side of (1.1.3) is isomorphic to
$$
\{\Gamma  \in  C^{1-n}(X, Y)_{{\Bbb Q}} \,| \,\Gamma \scirc
\Gamma _{\xi } = 0\,\,\, \text{and }{}^{t}\Gamma _{\xi '}\scirc
\Gamma  = 0\}.
\leqno(1.1.4)
$$

\noindent
{\bf 1.2.}
{\it Murre's construction.} Let
$ X $ be a connected smooth projective variety of dimension
$ n \ge 2 $.
We choose and fix an embedding of
$ X $ into a projective space.
Let
$ l $ denote the multiplication by the hyperplane section class.
Let
$ C $ be the intersection of
$ n - 1 $ generic smooth hyperplane sections.
(Note that
$ [C] \in  \CH^{n-1}(X) $ is independent of the choice of
$ C $.)
By (1.1.3--4) there exists uniquely
$ \Gamma  \in  \CH^{1}(X\times X)_{{\Bbb Q}} $ such that
$$
\align
\Gamma \scirc\Gamma _{\xi } &= 0,\quad {}^{t}\Gamma _{\xi '}\scirc
\Gamma  = 0,
\tag 1.2.1
\\
\Gamma _{*}\scirc l^{n-1} &= id\,\,\, \text{on }H^{1}(X,{\Bbb Q}).
\tag 1.2.2
\endalign
$$
We have
$ {}^{t}\Gamma  = \Gamma  $ if
$ \xi  = \xi ' $.
Note that (1.2.1) implies
$$
\Gamma _{*} : H^{i+2n-2}(X,{\Bbb Q})(n-1) \rightarrow  H^{i}(X,{\Bbb Q}
)\,\,\, \text{vanishes for }i \ne  1.
\leqno(1.2.3)
$$

Let
$ i : C \rightarrow  X $ denote the inclusion morphism.
Let
$$
\pi ' = \Gamma \scirc i_{*}\scirc i^{*},
$$
where
$ i_{*} = {}^{t}\Gamma _{i}, i^{*} = \Gamma _{i} $ as in (1.1.2).
Following [12] we define
$$
\pi _{0} = \Gamma _{[X]}\scirc {}^{t}\Gamma _{\xi '},\quad \pi _{1} = \pi 
'\scirc(1 - {}^{t}\pi '/2),\quad \pi _{2n-1} = {}^{t}\pi _{1},\quad \pi _{2n} = {}^{t}\pi 
_{0},
$$
where
$ 1 $ denotes the diagonal of
$ X $.
Then
$$
(\pi _{i})_{*}|H^{j}(X,{\Bbb Q}) = \delta _{i,j}id\quad  \text{for }
i = 0, 1, 2n-1, 2n.
\leqno(1.2.4)
$$
If
$ n = 2 $,
we define
$ \pi _{2} = 1 - \sum _{i\ne 2} \pi _{i} $.

\proclaim{{\bf 1.3.~Theorem}}
{\rm (Murre [7]).}
$ \pi _{i}\scirc \pi _{j} = \delta _{i,j}\pi _{i}\quad  $ for
$ \{i,j\} \subset  \{0,1, 2n-1, 2n\} $.
\endproclaim

\noindent
{\it Outline of proof.} We recall here some arguments of the proof which
will be needed in the proof of the main theorems.
See [7], [12] for details.

We have
$ \Gamma \scirc i_{*}\scirc i^{*}\scirc\Gamma  = \Gamma 
 $ by (1.1.3--4), and
$ {}^{t}\Gamma \scirc\Gamma  = 0 $ in
$ \CH^{2-n}(X\times X)_{{\Bbb Q}} $,
because it is cohomologically zero by (1.2.3).
So we get
$$
\pi ^{\prime 2} = \pi ',\quad {}^{t}\pi '\scirc\pi ' = 0.
\leqno(1.3.1)
$$
Then we can verify
$$
{\pi }_{1}^{2} = \pi _{1},\quad \pi _{2n-1}\scirc \pi _{1} = 0,\quad \pi _{1}\scirc \pi 
_{2n-1} = 0.
\leqno(1.3.2)
$$

We have furthermore
$$
\pi _{0}\scirc \pi ' = \pi _{2n}\scirc \pi ' = \pi '\scirc
\pi _{0} = \pi '\scirc\pi _{2n} = 0
\leqno(1.3.3)
$$
Indeed,
$ \pi _{0}\scirc \pi ' = 0  $ by (1.2.1), and
$ \pi '\scirc\pi _{0} = \pi _{2n}\scirc \pi ' = 0 $,
because
$ \pi '\scirc\Gamma _{[X]} \in  C^{0}(pt,X) $ is cohomologically zero (and
similarly for
$ {}^{t}\Gamma _{[X]}\scirc \pi ') $.
Finally the vanishing of
$ \pi '\scirc\pi _{2n} $ follows from
$ i_{*}\scirc i^{*}\scirc\Gamma _{\xi '} \in  C^{2n-1}(pt,
X) = 0 $.

Then we can verify the remaining assertions.

\medskip
\noindent
{\it Remark.} The Albanese map
$ \CH_{0}(X)_{{\Bbb Q}}^{0} \rightarrow  \Alb_{X}({\Bbb C})_{{\Bbb Q}} $ induces an isomorphism

$$
(\pi _{2n-1})_{*}\CH_{0}(X)_{{\Bbb Q}}^{0} \simto \Alb_{X}({\Bbb C})_{{\Bbb Q}}.
$$
If
$ n = \dim X = 2 $,
$ (\pi _{2})_{*}\CH_{0}(X)_{{\Bbb Q}}^{0} $ coincides with the kernel of the Albanese map with
$ {\Bbb Q} $-coefficients.
See [7, 7.1]. 

\medskip
\noindent
{\bf 1.4.}
{\it Construction of}
$ \tilde{\pi}_{2} $.
Assume
$ n = \dim X = 2 $ and
$ p_{g} = 0 $.
Let
$ C_{i} $ be irreducible curves on
$ X $ such that the cohomology classes of
$ [C_{i}] $ form a basis of
$ H^{2}(X,{\Bbb Q})(1) $.
Let
$ \tilde{C}_{i} $ be the normalization of
$ C_{i} $,
and
$ \tilde{C} $ the disjoint union of
$ \tilde{C}_{i} $ with
$ \tilde{i} : \tilde{C} \rightarrow  X $ the canonical morphism.
Let
$ A = (A_{i,j}) $ be the intersection matrix of the
$ C_{i} $ (i.e.
$ A_{i,j} = C_{i}\sbull C_{j}) $.
Let
$ B = (B_{i,j}) $ be the inverse of
$ A $.

Let
$ \Gamma _{B} \in  C^{-1}(\tilde{C},\tilde{C})_{{\Bbb Q}} = \CH^{0}(\tilde{C}
\times \tilde{C})_{{\Bbb Q}} $ defined by the matrix
$ B $.
Let
$$
\tilde{\Gamma} = \tilde{i}_{*}\scirc \Gamma _{B}\scirc \tilde{i}^{*}.
$$
Since the composition
$$
\tilde{i}^{*}\scirc\tilde{i}_{*} : H^{0}(\tilde{C},{\Bbb Q}) \rightarrow  
H^{2}(\tilde{C},{\Bbb Q})(1)
$$
is given by the matrix
$ A $ (using the projection formula),
we see that
$$
\tilde{\Gamma}_{*}|H^{i}(X,{\Bbb Q}) = \delta _{i,2}id. 
\leqno(1.4.1)
$$
(The assertion for
$ i \ne  2 $ follows from the definition of
$ \tilde{\Gamma} $ .)
Note that
$ {}^{t}\tilde{\Gamma} = \tilde{\Gamma} $ because
$ B $ is symmetric.
We define
$$
\tilde{\pi}_{2} = (1 - \pi ')\scirc\tilde{\Gamma}\scirc
(1 - {}^{t}\pi ')
$$
Then the symmetry
$ {}^{t}\tilde{\pi}_{2} = \tilde{\pi}_{2} $ is clear.

\proclaim{{\bf 1.5.~Proposition}}
$ \tilde{\pi}_{2}  $ is an idempotent, and is orthogonal to
$ \pi _{i} $ for
$ i \ne  2 $.
\endproclaim

\demo\nofrills {Proof.\usualspace}
We have
$ \pi '\scirc\tilde{\pi}_{2} = 0 $ by (1.3.1).
Since
$$
\Gamma _{B}\scirc \tilde{i}^{*}\scirc\Gamma  \in  C^{-2}(X,
\tilde{C})_{{\Bbb Q}} = \CH^{0}(X\times \tilde{C})_{{\Bbb Q}}
$$
is cohomologically zero by (1.2.3), we get
$$
\tilde{\Gamma}\scirc\pi ' = 0,\quad \tilde{\pi}_{2}\scirc \pi ' = 0,
\leqno(1.5.1)
$$
using (1.3.1).
Then we have
$ {}^{t}\pi '\scirc\tilde{\pi}_{2} = \tilde{\pi}_{2}\scirc {}^{t}\pi '=0 $ by 
transpose, and
$$
\pi _{i}\scirc \tilde{\pi}_{2} = \tilde{\pi}_{2}\scirc \pi _{i} = 0\quad 
 \text{for }i = 1, 3.
\leqno(1.5.2)
$$

For
$ i = 0, 4 $,
we have
$ \tilde{\Gamma}\scirc\pi _{0} = \tilde{\Gamma}\scirc\pi _{4} 
=0 $ by
$ \Gamma _{B}\scirc \tilde{i}^{*}\scirc\Gamma _{[X]} \in  C
^{-1}(pt,\tilde{C}) = 0 $ and
$ \tilde{i}^{*}\scirc\Gamma _{\xi '} \in  C^{2}(pt,\tilde{C}) = 0 $.
Then (1.5.2) holds also for
$ i = 0, 4 $ using (1.3.3), because
$ {}^{t}\tilde{\Gamma} = \tilde{\Gamma} $.

Finally we have
$$
{\tilde{\pi}}_{1}^{2} = \tilde{\pi}_{2}.
\leqno(1.5.3)
$$
Indeed,
$ \tilde{\Gamma} $ is an idempotent, because
$ \Gamma _{B}\scirc\tilde{i}^{*}\scirc\tilde{i}_{*}\scirc
\Gamma _{B} \in  \CH^{0}(\tilde{C}\times \tilde{C})_{{\Bbb Q}} $ 
coincides with
$ \Gamma _{B} $ (cohomologically).
Then (1.5.3) follows from (1.5.1) and (1.3.1).
\enddemo

\noindent
{\bf 1.6.}
{\it Proof of} (0.1) 
$ \Rightarrow  (0.3) $.
Applying (0.1) to 
$ 1 - (\sum_{i \ne 2} \pi _{i} + \tilde{\pi}_{2}) $, we get
$ \tilde{\pi}_{2} = \pi _{2} $ by (0.1).
So it is enough to show
$ (\tilde{\pi}_{2})_{*}\CH^{2}(X)_{{\Bbb Q}} = 0 $ by [7, 7.1] (because the Albanese kernel is torsion free by [9]).
See Remark after (1.3).
By the definition of
$ \tilde{\pi}_{2} $,
it is enough to show
$ \tilde{i}^{*}\CH^{2}(X)_{{\Bbb Q}} = 0 $.
But this is clear.

\bigskip\bigskip\centerline{{\bf 2. Cycle maps and correspondences}}

\bigskip
\noindent
{\bf 2.1.}
Let
$ X $ be a smooth proper complex variety with the structure morphism
$ a_{X} : X \rightarrow  pt $.
Let
$ {\Bbb Q} (j) $ denote the Tate Hodge structure of type
$ (-j,-j) $ (see [5]) which is naturally identified with a mixed Hodge Module
on
$ pt  $ (see [10, (4.2.12)]).
Then we have a cycle map
$$
cl : \CH^{p}(X)_{{\Bbb Q}} \rightarrow  \Ext^{2p}({a}_{X}^{*}{\Bbb Q}, {a}_{X}^{*}
{\Bbb Q}(p)) = \Ext^{2p}({\Bbb Q}, (a_{X})_{*}{a}_{X}^{*}{\Bbb Q}(p)),
\leqno(2.1.1)
$$
where
$ {\Bbb Q} $ means
$ {\Bbb Q}(0) $, and
$ \Ext $ is taken in the derived category of mixed Hodge Modules on
$ X $ or
$ pt $.
See (4.5.18) of loc. cit.  
(The last isomorphism of (2.1.1) follows from the adjoint relation.)
The target is isomorphic to
$ {\Bbb Q} $-Deligne cohomology, and (2.1.1) is an isomorphism for
$ p = 1 $ as well-known.
(This cycle map coincides with Deligne's cycle map which uses local 
cohomology, and can also be obtained by using the theory of Bloch and 
Ogus as was done by Beilinson and Gillet.
In particular, its restriction to homologically equivalent to zero cycles 
coincides with Griffiths' Abel-Jacobi map.)

Let
$ X, Y $ be smooth proper complex varieties.
Assume
$ X $ is purely
$ n $-dimensional.
Then the cycle map induces
$$
\aligned
cl : C^{i}(X,Y)_{{\Bbb Q}} &\rightarrow  \Ext^{2n+2i}({\Bbb Q}, (a
_{X\times Y})_{*}{a}_{X\times Y}^{*}{\Bbb Q}(n+i))
\\
&=\Ext^{2i}((a_{X})_{*}{a}_{X}^{*}{\Bbb Q}, (a_{Y})_{*}{a}_{Y}^{*}{\Bbb Q} \text{(i)) } .
\endaligned
\leqno(2.1.2)
$$
See [11, II].
This is an isomorphism for
$ n + i = 1 $.
By (3.3) of loc. cit, (2.1.2) is compatible with the composition of correspondences.

\proclaim{{\bf 2.2.~Proposition}}
Let 
$ X, S $ be connected smooth proper complex varieties, and
$ \Gamma  \in  \CH^{p}(S\times X)_{{\Bbb Q}} $.
Assume
$ \Gamma  $ is homologically equivalent to zero, and the cycle map (2.1.1) 
for
$ X $ and
$ p $ is injective.
Then
$ \Gamma  = \Gamma _{1} + \Gamma _{2} $ where
$ \Gamma _{1} $ is supported on
$ D \times  X $ for a divisor
$ D $ on
$ S $,
and
$ \Gamma _{2} = [S] \times  \xi  $ for
$ \xi  \in  \CH^{p}(X)_{{\Bbb Q}} $.
\endproclaim

\demo\nofrills {Proof.\usualspace}
Let
$ H^{2p-1}(X,{\Bbb Q})_{S} $ denote a constant variation of Hodge structure
on
$ S $ such that the fibers are
$ H^{2p-1}(X,{\Bbb Q}) $.
This is identified with the direct image of
$ {\Bbb Q}_{S\times X} $ by the first projection.
Since the restriction of
$ \Gamma  $ to
$ \{s\} \times  X $ is homologically equivalent to zero for any
$ s \in  S $,
$ \Gamma  $ determines a normal function
$$
e \in  \Ext^{1}({\Bbb Q}_{S}, H^{2p-1}(X,{\Bbb Q})_{S}(p)),
$$
where
$ \Ext $ is taken in the derived category of mixed Hodge Modules (using [10, 
3.27]).
Note that
$ e $ can be identified with a section of
$ S\times J^{p}(X) \rightarrow  S $ (where
$ J^{p}(X) $ is Griffith' intermediate Jacobian), if we replace rational
coefficients with integral coefficients.

By the adjoint relation for
$ a_{S} : S \rightarrow  pt $,
we have a short exact sequence
$$
\aligned
0 \rightarrow  \Ext^{1}({\Bbb Q}, H^{2p-1}(X,{\Bbb Q})
&(p)) \rightarrow  \Ext^{1}({\Bbb Q}_{S}, H^{2p-1}(X,{\Bbb Q})_{S}(p))
\\
&\rightarrow  \Hom({\Bbb Q}, H^{1}(S,{\Bbb Q})\otimes H^{2p-1}
(X,{\Bbb Q})(p)) \rightarrow  0,
\endaligned
$$
where the first
$ \Ext $ and the last
$ \Hom  $ are taken in the category of mixed Hodge structures.
Since
$ \Gamma  $ is homologically equivalent to zero, the image of
$ e $ in the last term is zero, and hence
$ e $ comes from the first term, i.e.
it is constant.
So there exists
$ \xi  \in  \CH^{p}(X)_{{\Bbb Q}} $ such that
$ \xi  $ is homologically equivalent to zero and
$ e $ is the image of
$ [S] \times  \xi  $.
Then replacing
$ \Gamma  $ with
$ \Gamma  - [S] \times  \xi  $,
we may assume
$ e = 0 $.

Let
$ k $ be a finitely generated subfield of
$ {\Bbb C} $ such that
$ X, S $ and
$ \Gamma  $ are defined over
$ k $,
i.e.
there exist smooth proper
$ k $-varieties
$ X_{k}, S_{k} $ and
$ \Gamma _{k} \in  \CH^{p}(S_{k} \times _{k} X_{k})_{{\Bbb Q}} $ with isomorphisms
$ X = X_{k} \otimes _{k} {\Bbb C},  $ etc.
Let
$ K = k(S_{k}) $ the function field of
$ S_{k} $,
and
$ X_{K} $ the generic fiber of the first projection of
$ S_{k} \times _{k} X_{k} $.
Let
$ \Gamma _{K} \in  \CH^{p}(X_{K})_{{\Bbb Q}} $ denote the restriction of
$ \Gamma _{k} $.
Then it is enough to show
$ \Gamma _{K} = 0 $.

We choose an embedding
$ K \rightarrow  {\Bbb C} $ extending
$ k \rightarrow  {\Bbb C} $.
Since
$ e = 0 $,
the image of
$ \Gamma _{K}\otimes _{K}{\Bbb C} \in  \CH^{p}(X)_{{\Bbb Q}}
 $ by the cycle map is zero, because
$ \Gamma _{K}\otimes _{K}{\Bbb C} $ is identified with the restriction
of
$ \Gamma  $ to
$ \{s\}\times X $ where
$ s \in  S $ corresponds to the embedding
$ K \rightarrow  {\Bbb C} $.
So
$ \Gamma _{K}\otimes _{K}{\Bbb C} $ is zero by hypothesis.
Then the assertion follows from the injectivity of
$ \CH^{p}(X_{K})_{{\Bbb Q}} \rightarrow  \CH^{p}(X)_{{\Bbb Q}} $.
\enddemo

\proclaim{{\bf 2.3.~Theorem}}
Let 
$ X $ be a connected smooth proper complex surface such that the Albanese map for
$ X $ is injective.
Then the cube of the ideal in (0.2) is zero.
\endproclaim

\demo\nofrills {Proof.\usualspace}
Let
$ \Gamma  \in  \CH^{2}(X \times  X)_{{\Bbb Q}} $ that is homologically equivalent
to zero.
Then
$ \Gamma  = \Gamma _{1} + \Gamma _{2} $ such that
$ \Gamma _{1} $ is supported on
$ D \times  X  $ and
$ \Gamma _{2} = [X] \times  \xi  $ as in (2.2).
We have to show
$$
\Gamma ''\scirc\Gamma '\scirc\Gamma _{i} = 0\,\,\,(i = 1, 2)
\leqno(2.3.1)
$$
for any
$ \Gamma ', \Gamma''  \in  C^{0}(X,X) $ which are homologically equivalent
to zero.

For
$ i = 1 $,
let
$ Y $ denote the normalization of
$ D $ with
$ f : Y \rightarrow  X $ the canonical morphism.
Then there exists
$ \Gamma '_{1} \in  C^{0}(Y,X)_{{\Bbb Q}} $ such that
$ \Gamma _{1} = \Gamma '_{1}\scirc f^{*} $.
By the injectivity of the cycle map in the divisor case, it is enough to show
that the image of
$ \Gamma ''\scirc\Gamma '\scirc\Gamma '_{1} $ by (2.1.2) is zero.
Since
$ \Gamma ', \Gamma '' $ are homologically equivalent to zero, and
$ \Ext^{2} $ vanishes, the assertion follows by using for example a (noncanonical) decomposition
$$
(a_{X})_{*}{a}_{X}^{*}{\Bbb Q} \simeq \oplus _{i} H^{i}(X,{\Bbb Q})[-i]
$$
in the derived category of mixed Hodge Modules on
$ pt $ (or equivalently, of graded-polarizable mixed Hodge structures).
See [10, (4.5.4)].

The argument is similar for
$ i = 2 $ by using the injectivity of the Albanese map, because it is enough
to show that the image of
$ \Gamma ''\scirc \Gamma '\scirc \Gamma _{\xi } \in  
C^{2}(pt,X) $ by the cycle map (2.1.2) is zero.
This completes the proof of (2.3).
\enddemo

\bigskip
\noindent
{\bf 2.4.}
{\it Remark.} The square of the ideal in (0.2) is nonzero if
$ H^{1}(X,{\Bbb Q}) \ne  0 $.
Indeed, let
$ C $ be a hyperplane section of
$ X $ with the inclusion morphism
$ i : C \rightarrow  X $.
By Hodge theory we have a divisor
$ D $ on
$ C\times X $ such that
$ D_{*^{}} : H^{i}(C,{\Bbb Q}) \rightarrow  H^{i}(X,{\Bbb Q}) $ is zero for
$ i \ne  1 $ and
$ D_{*}\scirc  i^{*} : H^{1}(X,{\Bbb Q}) \rightarrow  H^{1}(X,{\Bbb Q}
) $ is an isomorphism.
Let
$ k $ be a finitely generated subfield of
$ {\Bbb C} $ such that
$ X $,
$ C $,
$ D $ are defined over
$ k $,
i.e.
there exist
$ X_{k} $,
$ C_{k} $,
$ D_{k} $ with isomorphisms
$ X_{k}\otimes _{k}{\Bbb C} = X,  $ etc.
Let
$ D_{i} = (pr_{i}\times id)^{*}D_{k} $ where
$ pr_{i} : C_{k}\times _{k}C_{k} \rightarrow  C_{k} $ are the natural projections
for
$ i = 1, 2 $.
Let
$ K = k(C_{k}\times _{k}C_{k}) $,
and
$ \Gamma _{K,i} $ be the restriction of
$ D_{i} $ to the generic fiber of
$ C_{k}\times _{k}C_{k}\times _{k}X_{k} \rightarrow  C_{k}\times _{k}C
_{k} $.
Then
$ \Gamma _{K,1}\times _{K}\Gamma _{K,2} $ is the restriction of
$ D_{k}\times _{k}D_{k} $.
Let
$ \Gamma _{i} = \Gamma _{K,i}\otimes _{K}{\Bbb C} \in  C^{1}(pt
,X)_{{\Bbb Q}} $.
Then the composition of
$ \Gamma _{[C]}\scirc {}^{t}\Gamma _{1} $ and
$ \Gamma _{1}\scirc {}^{t}\Gamma _{[C]} $ is equal to a nonzero multiple
of
$ \Gamma _{2}\scirc {}^{t}\Gamma _{1} $,
and this is nonzero.

\bigskip\bigskip\centerline{{\bf References}}

\bigskip

\item{[1]}
R. Barlow, Rational equivalence of zero cycles for some more surfaces
with
$ p_{g} = 0 $, Inv.  Math. 79 (1985), 303--308.

\item{[2]}
A. Beilinson,
Height pairing between algebraic cycles, Lect. Notes in Math., vol.
 1289, Springer, Berlin, 1987, pp. 1--26.

\item{[3]}
S. Bloch, Lectures on algebraic cycles, Duke University Mathematical
series 4, Durham, 1980.

\item{[4]}
S. Bloch, A. Kas and D. Lieberman, Zero cycles on surfaces with
$ p_{g} = 0,  $ Compos. Math. 33 (1976), 135--145.

\item{[5]}
P. Deligne, Th\'eorie de Hodge I, Actes Congr\`es Intern. Math.,
1970, vol. 1, 425-430 : II,  Publ. Math. IHES, 40 (1971), 5--57; III ibid., 
44 (1974), 5--77.

\item{[6]}
D. Mumford, Rational equivalence of
$ 0 $-cycles on surfaces, J. Math. Kyoto Univ. 9 (1969), 195--204.

\item{[7]}
J.P. Murre, On the motive of an algebraic surface, J. Reine Angew. Math. 409 (1990), 
190--204.

\item{[8]}
\SameAuthor, On a conjectural filtration on Chow groups of an algebraic variety,
 Indag. Math. 4 (1993), 177--201.

\item{[9]}
A. Roitman, The torsion in the group of zero cycles modulo rational equivalence,
 Ann. Math. 111 (1980), 553--569.

\item{[10]}
M. Saito, Mixed Hodge Modules, Publ. RIMS, Kyoto Univ., 26 (1990), 
221--333.

\item{[11]}
\SameAuthor, Hodge conjecture and mixed motives, I, Proc. Symp. Pure Math.
53 (1991), 283--303; II, in Lect. Notes in Math., vol. 1479, Springer, Berlin, 1991, pp.
 196--215.

\item{[12]}
A. J. Scholl, Classical Motives, Proc. Symp. Pure Math.
55 (1994), Part 1, 163--187.

\item{[13]}
C. Voisin, Sur les z\'ero cycles de certaines hypersurfaces munies
d'un automorphisme, Ann. Sci. Norm. Sup. Pisa 19 (4) (1992), 473--492.

\bigskip
\noindent
Feb. 7, 2000
\bye